\documentclass[a4paper,12pt]{amsart}
\usepackage{bm}
\usepackage{amsfonts,latexsym,rawfonts,amsmath,amssymb,amsthm, mathrsfs, lscape,color}
\usepackage{extsizes}
\usepackage[all]{xy}
\usepackage[english]{babel}		\usepackage[T1]{fontenc}
\usepackage{graphicx}
\usepackage{booktabs}
\usepackage{array, tabularx}
\usepackage{setspace}
\usepackage{textcomp}

\usepackage{tikz}
\usepackage{multirow}
\usepackage{paralist}
\usepackage{comment}

\usepackage[hypertexnames=false,
backref=page,
    pdftex,
    pdfpagemode=UseNone,
    breaklinks=true,
    extension=pdf,
    colorlinks=true,
    linkcolor=blue,
    citecolor=blue,
    urlcolor=blue,
]{hyperref}

\usepackage[top=1.5in, bottom=.9in, left=0.9in, right=0.9in]{geometry}

\newcommand{\referenza}{}

\newtheorem{thm}{Theorem}[section]
\newtheorem*{thm*}{Theorem \referenza}
\newtheorem{cor}[thm]{Corollary}
\newtheorem*{cor*}{Corollary \referenza}
\newtheorem{lem}[thm]{Lemma}
\newtheorem*{lem*}{Lemma \referenza}

\newtheorem{prop}[thm]{Proposition}
\newtheorem*{prop*}{Proposition \referenza}
\newtheorem{prob}[thm]{Problem}

\newtheorem*{conj*}{Conjecture \referenza}
\newtheorem{rmk}[thm]{Remark}
\newtheorem*{rmk*}{Remark \referenza}

\newtheorem{exa}[thm]{Example}
\newtheorem{defi}[thm]{Definition}
\newtheorem{exe}[thm]{Exercise}

\def\bexc{\begin{exe}}
\def\eexc{\end{exe}}

\def\be{\begin{equation}}
\def\ee{\end{equation}}
\def\bp*{\begin{prop*}}
\def\ep*{\end{prop*}}
\def\bt*{\begin{thm*}}
\def\et*{\end{thm*}}
\def\bt{\begin{thm}}
\def\et{\end{thm}}
\def\bp*{\begin{prop*}}
\def\ep*{\end{prop*}}
\def\bp{\begin{prop}}
\def\ep{\end{prop}}
\def\bl*{\begin{lem*}}
\def\el*{\end{lem*}}
\def\bl{\begin{lem}}
\def\el{\end{lem}}
\def\bec*{\begin{cor*}}
\def\enc*{\end{cor*}}
\def\bc{\begin{cor}}
\def\ec{\end{cor}}
\def\bex{\begin{exa}}
\def\eex{\end{exa}}
\def\de{\begin{defi}}
\def\ed{\end{defi}}
\def\br{\begin{rmk}}
\def\er{\end{rmk}}
\def\pb{\begin{prob}}
\def\eb{\end{prob}}

\def\ES{\varnothing}
\numberwithin{equation}{section}

 \def \t\G {\widetilde \Gamma} \def \t\T {\widetilde \T}
\def\bit{\begin{itemize}} \def\eit{\end{itemize}}
\def\benu{\begin{enumerate}} \def\enu{\end{enumerate}}
\def\demo{\begin{proof}} \def\enddemo{\end{proof}}

\def \p {\partial}
\def\de{\p}

\numberwithin{equation}{section}
\def\acal{\mathcal A}    
\def\hcal{\mathcal H}\def\lcal{\mathcal L}\def\mcal{\mathcal M}
\def\ocal{\mathcal O}\def\pcal{\mathcal P}

\def\lsf{{\sf L}}

\def\R{{\mathbb {R}}}   \def\a {\alpha} \def\b {\beta}
\def\N{{\mathbb N}}     \def\d {\delta} \def\e{\varepsilon}
\def\C{{\mathbb C}}      
					\def\p{\partial}
 \def\O{\Omega}
					\def\t{\theta}\def\z{\zeta}

\def\T{{\mathbb T}}

 \def\G{\Gamma}
\def\O{\Omega}
\def\R{{\mathbb {R}}}   \def\a {\alpha} \def\b {\beta}
 
\def\N{{\mathbb N}}     \def\d {\delta} \def\e{\varepsilon}

\def\B{{\mathbb B}}
\def\C{{\mathbb C}}

\def\G{{\mathbb G}}

\def\R{{\mathbb R}}

					\def\p{\partial}

					\def\t{\theta}\def\z{\zeta}

\def\T{{\mathbb T}}

   \def\a {\alpha} \def\b {\beta}
\def\N{{\mathbb N}}     \def\d {\delta} \def\e{\varepsilon}
\def\C{{\mathbb C}}      
				\def\p{\partial}

				\def\t{\theta}\def\z{\zeta}

 \def\G{\Gamma}

\def\gl{\gtrless}
\def\1{\rm 1} \def\2{\rm 2} \def\3{\rm 3} \def\4{\rm 4} \def\5{\rm 5} \def\6{\rm 6}

\def\lp{\rm(\,} \def\rp{\,\rm)}
\def\smi{\setminus} \def\ssmi{\!\smallsetminus\!}
\def\til{\tilde}

\def\Til {\widetilde}
\def \Hat {\widehat}

\def\IN{\infty}

\def\tms{\times}

\def\Int{\stackrel{\rm o}}
 \def\oli{\overline}
\def\sbs{\subset} 
\def\sbseq{\subseteq}  
\def\nsbs{\nsubseteq}  
\def\Sbs{\Subset} 
 
\def\beqn{\begin{eqnarray*}} %non numerata
\def\eeqn{\end{eqnarray*}}
\def\beqnn{\begin{eqnarray}} % numerata
\def\eeqnn{\end{eqnarray}}
\def\ba{\begin{aligned}}
\def\ea{\end{aligned}}
\def\bca{\begin{cases}}
\def\eca{\end{cases}}
 
\def\neqv{\not\equiv}
 
\def\nin{\noindent}
\def\bit{\begin{itemize}}
\def\eit{\end{itemize}}

\allowdisplaybreaks[1]

\title[Levi equation and local maximum property]
{Levi equation and local maximum property}

\author{Giuseppe Della Sala and Giuseppe Tomassini}
\address[]{Department of Mathematics, Faculty of Arts and Sciences, AUB Riad El Solh, Beirut 1107 2020, Lebanon
}
\email{gd16@aub.edu.lb} 
\address[]{Scuola Normale Superiore, Piazza dei Cavalieri, 7 - I-56126 Pisa, Italy
}              
\email{giuseppe.tomassini@sns.it}
\thanks{}

\date{\today}

\begin{document}
\begin{abstract} The aim of the paper is to study the level sets of the solutions of Dirichlet problems for the Levi operator on strongly pseudoconvex domains $\Omega$ in $\mathbb C^2$. Such solutions are generically non smooth, and the geometric properties of their level sets are characterized by means of hulls of their intersections with $b\Omega$, using as main tool the local maximum property introduced in \cite{S1}. The same techniques are then employed to study the behavior of the complete Levi operator for graphs in $\mathbb C^2$.
 \end{abstract}

\maketitle
\tableofcontents
\section*{Introduction}\label{INTR}

The principal aim of the paper is to study the relationship between the level sets of weak solutions for the \emph{Levi equation} 
\beqnn\label{5(6)} 
&&{\lcal
}(u):=-{\rm \det}\left(
\begin{array}{ccc}
0&u_{z_1}&u_{z_2}\\ u_{\oli z_1}&u_{\oli z_1 z_1}&u_{\oli z_1 z_2}\\ u_{\oli z_2}&u_{\oli z_2 z_1}&u_{\oli z_2 z_2} 
\end{array}\right)=\\
&&\sum_{\a,\b=1}^2(\delta_{\a\b}{\vert\partial
u\vert^{2}-u_{\bar\a}u_{\b}})u_{\a\bar\b}=0\nonumber
\eeqnn
($u_\a=u_{z\a}$,  $u_{\bar\a}=u_{{\bar z}_\a}$)
on a strongly pseudoconvex domain $\Omega\subset \C^2$. It is well known that the regular levels of a smooth solution of $\lcal(u)=0$ are Levi flat  (i.e. foliated by complex curves). This is no longer true for (continuous weak)  solutions in the sense of viscosity \cite{CIL}. However, as discovered in \cite{ST}, the zero sets of these solutions retain an interesting property, the ``local maximum property''. The subsets of $\C^2$ (or $\C^n$) endowed with this property form a category which contains properly the complex analytic and Levi flat subsets. They were introduced in Complex Analysis by Slodkowski (see  \cite{S1})  who demonstrated two basic results that have revealed fundamental tools for our research: a version of the classical ``Kontinuit\"atsatz''  and a theorem of duality between pseudoconvexity and local maximum property (see Sec. 2). 

The paper must be considered as a natural continuation  of \cite{ST}. It is organized in three sections. In Sec. 1 we summarize the main results of \cite{ST}, particularly those related to the Dirichlet problem
\[
\hspace{-1cm}({\sf P}_{g})\hspace{1cm}\begin{cases}
\lcal(u)=0 \>\> & {\rm in}\>\> \O\\
u=g\>\> & {\rm on}\>\> {\rm b}\O
\end{cases}
\] 
with $g$ continuous. This problem has a continuous weak solution  $u$ which is not unique in general. Neverthless all the solutions share the following properties:
\bit
\item[1\rp ] every level set $\{u=c\}$ is l.m.p.
\item[2\rp ] $\{u=c\}\sbseq{\Hat{\{g=c\}}}_\hcal$ 
\eit
where $\hcal=\hcal(\O)$ is the Banach algebra $\rm{C}^0\cap\ocal(\O)$ and ${\Hat{\{g=c\}}}_\hcal$  denotes the $\hcal$-hull of $\{g=c\}$ (see Sec. 2).

The geometric counterpart of this theorem tells us that every closed subset $S$ of ${\rm b}\O$ is the boundary of a compact subset $\oli M$ of $\oli\O$ such that $M:=\oli M\cap\O$ is an l.m.p subset. In particular, $M$ is a member of a continuous family of l.m.p. subsets of $\O$. 

In Section 2 we focus our attention on the inverse problem which is the bulk of the present paper : given a pair $(S,\oli M)$ as above put it into a continuous family $\{\oli M_t\}$ of closed subsets of $\oli\O$ such that 
$\oli M\sbseq \oli M_0$ and the sets $\oli M_t\cap\O$ are l.m.p. The main results  are contained in the theorems \ref{T10}, \ref{C21} obtained studying an appropriate Dirichlet problem $({\sf P}_{g})$. Precisely, we assume that  ${\rm }b\O\ssmi S$ has two connected components  $S^\pm$ and we can suppose $S=\{g= 0\}$ where $g$ is a continuous function $g:{\rm b}\O\to\R$ such that $S^\pm=\{g\gl 0\}$ and $g$ is {\em minimal} i.e. $g$ has only two peak points $p\in S^+$, $q\in S^-$ (see Sec. \ref{MT}). Then if $u\in{\rm C^0}(\oli\O)$ is a weak solution of the problem $({\sf P}_{g})$ we prove (see Theorem \ref{T10}) that
$$ 
\rm{a)}\>\>\> M\sbseq\{u=0\}\sbseq \Hat{\{g=0\}}_{\hcal}=\Hat{S}_{\hcal}. 
$$
Moreover, if $S$ is ${\rm C}^1$ and $\O\ssmi M$ has two connected components $\O^\pm$, then
\beqn
b)&&M\sbseq\{u=0\}\sbseq\Hat M_{\oli\ocal^+}\cup\Hat M_{\oli\ocal^-}=\\
&&{\Hat S}_{\oli\ocal^+}\cup\Hat {S}_{\oli\ocal^-}\sbseq\Hat{S}_{\oli\ocal}
\eeqn  
where $\oli{\ocal}^\pm=\ocal(\oli{\O^\pm})$, $\oli\ocal=\ocal(\oli\O).$
\begin{rmk*}
From $\rm{a)}$ it follows in particular that there is a weak solution $u=u_g\neqv 0$ of the Levi equation such that $u=0$ on $M$. We will compare this result with the main result of \cite{FO}.
\end{rmk*}
We can rephrase the above by the following statement (see Theorem \ref{C21})
\begin{thm*}
Let $S$ be a compact subset of ${\rm b}\O$ and $\mcal(S)$ the family of all closed subsets $M$ of $\oli\O$ such that $M\cap\O$ is l.m.p. and $M\cap{\rm b}\O=S$. If $S$ divides ${\rm b}\O$ into two connected components then
\bit
\item[1\rp]$\mcal(S)$ has a unique maximal element $M_S;$
\item[2\rp] $M_S=\{u_g=0\}$ where $g$ is a minimal defining function for $S$ and $u_g$ is an arbitrary solution of  solution of $({\sf P}_{g})$;
\item[3\rp]$M_S=\Hat S_\hcal.$
\eit 
\end{thm*}
\nin To establish part  3) we use that $\Hat S_\hcal$ is an l.m.p. subset (see Theorem \ref{hulmp}). For the proof of this fact we are indebited to Slodkowski.

From what is preceding it follows that if $M=\Hat{S}_{\pcal}$ then $M$ is the zero set of a weak solution $\{u=0\}$  of $({\sf P}_{g})$. In that case $M$ belongs to a continuous family of l.m.p. sets, and admits a Stein basis (see Corollary \ref{T11}):
\begin{thm*}
 There exists a basis of neighborhoods $\{U_n\}_{n\in\N}$ of $M$, for the relative topology of $\oli\O$ such that $U_n\cap\O$ and $\O\ssmi{\rm b}U_n$, $n\in\N$  are Stein domains
\end{thm*} 
The methods employed in the previous sections  apply in Section 3 to study properties of solutions of the  ``complete Levi operator''  $\lsf(v,k)$ for graphs (see  (\ref{L8}))  depending on the sign and the behaviour of $k$ .  Theorem \ref{T12} and Corollary \ref{T38}) are examples of the kind of results we can expect.
\section{Preliminaries}\label{Geo}
In this section we fix the notations and we recall the main results of \cite{ST}. 

Let $\O$ be a bounded domain in $\C^2$. We consider the function spaces
 $
 \ocal(\oli\O), \pcal(\O):={\rm Psh}(\O)\cap {\rm C^0(\oli \O})$ and  $\hcal(\O):=\ocal(\O)\cap {\rm C^0(\oli \O}). $
For a compact set $K\sbs\oli\O$ we denote $\Hat K_{\oli\ocal}$, $\Hat K_{\hcal}$ and $\Hat K_{\pcal}$ the respective hulls.  If $\O$ is strongly pseudoconvex,  for every compact subset $K$ of $\O$ the envelopes $\Hat K_{\pcal}$ and $ \Hat K_{\hcal}$ coincide (\cite[Theorem 1.2]{GS}).

Let us recall the definition of subsets with the local maximum property introduced in \cite{S1}.

A locally closed subset $X$ of $\C^n$, $n\ge 2$ has the {\em local maximum property} (or  $X$ is an {\em l.m.p.} subset) if for every point $x\in X$ there is a neighborhood $U$ of $x$ in $\C^n$ with following property: for every compact set $K\sbs U$ and every function $\phi$ which is plurisubharmonic in a neighborhood of $K$
\be\label{LM}
\max\limits_{K\cap X}\phi=\max\limits_{X\cap{\rm b}K}\phi
\ee
where $\max\limits_{X\cap{\rm b}K} \phi$ is meant to be $-\IN$ if $X\cap{\rm b}K\!\!=\ES$.

Accordingly, an  l.m.p. subset has no isolated point. 

It is immediate to check that if $V$ is a relatively open subset of $X$ with $\oli V$ compact and $\phi$ is plurisubharmonic in a neighborhood of $\oli V$ then 
$$
\max\limits_{\oli V}\phi=\max\limits_{{\rm b}\oli V}\phi.
$$
In \cite[Prop. 4.2]{S1} is proved that $X$ has the l.m.p. if and only if one of the two following equivalent conditions holds true:
\bit
\item[i)] condition \ref{LM} holds for every function $\vert P(z)\vert$ where $P(z)$ is a polynomial;
\eit
\bit
\item[ii)] there do not exist $p\in X$, $\e>0$ and a strongly plurisubharmonic function $\psi$ on the ball $B(p,r)=\{\vert z-p\vert<r\}$ such that $\psi(p)=0$ and $\psi(z)\le-\e\vert z-p\vert^2$ for $z\in X\cap B(p,r)$. 
\eit 
\subsection{Local maximum property and ``Kontinuit\"atsatz''}
Two facts are crucial for l.m.p subsets. The first one is a version of the classic ``Kontinuit\"atsatz'', proved  by Slodkowski in (\cite{S2}). Let $U$ be a domain in $\C^n$, $n\ge 2$, and $\{X_\nu\}_\nu$ a sequence of relatively compact subsets of $U$. Assume that every $X_\nu$ is l.m.p. and that
\bit
\item[iii)] $\limsup\limits_{\nu\to+\IN}(\oli {X_\nu}\ssmi X_\nu)$ is a compact subset of $U$
\item[iv)] there exists a sequence $\{x_\nu\}_\nu, x_\nu\in X_\nu$, converging to a point $x\in\C^n\ssmi U $;
\eit
then $U$ is not pseudoconvex.

The second one holds in $\C^2$ and concerns the duality between pseudoconvexity and l.m.p.. A relatively closed subset $X$ of a domain $V\sbs\C^2$ is l.m.p. if and only if $U=V\ssmi X$ is relatively pseudoconvex in $V$ (i.e. $U$ is locally pseudoconvex at the points of ${\rm b}U\cap V$). In particular, if $V$ is pseudoconvex then $X$ is l.m.p if and only if $U$ is pseudoconvex (\cite[Theorem 2]{S1}.
\subsection{The Dirichlet problem $\lp{\sf P}_g\rp$}\label{Exs} 
Here we summarize the main results  obtained in \cite{ST} for the Dirichlet problem
\[
\hspace{-1cm}({\sf P}_{g})\hspace{1cm}\begin{cases}
\lcal(u)=0 \>\> & {\rm in}\>\> \O\\
u=g\>\> & {\rm on}\>\> {\rm b}\O.
\end{cases}
\]
where $\O$ is a strongly pseudoconvex (or, more generally, $P$-regular, see \cite[Sec. 4]{ST})
 bounded domain in $\C^2$ and $g:{\rm b}\O\to\R$ is a continuous function. 

\bit
\item[1\rp]The  problem $\lp{\sf P}_g\rp$ has two extremal (continuous weak) solutions $u^\pm\in {\rm C}^0(\oli \O)$. All solutions of the same problem satisfy $u^-\le u \le u^+$ (\cite[Th. 3.3] {ST}).
\eit
\medskip
For a function $u\in {\rm C^0}( \O)$,
\bit
\item[2\rp] $u$ is a weak solution of ${\lcal}(u)=0$;
\item[3\rp] for every $c\in\R$, sets $\{u<c\}$, ($\{u\le c\}$), $\{u>c\}$, ($\{u\ge c\}$), ($\{u=c\}$) are pseudoconvex (l.m.p) whenever non-empty$\rp$ 
\eit 
are equivalent properties (\cite[Th. 3.3]{ST}).

Moreover, if $u\in {\rm usc}(\oli\O)$ is an upper semicontinuous  solution of ${\lcal
}(u)=0$ then 
\be\label{max}
\max\limits_{\oli\O}u=\max\limits_{\rm b\O}u
\ee
(\cite[Cor. 3.3]{ST}). In particular, if $u=0$ on ${\rm b}\O$ then $u\equiv 0$.
\br\label{R1}
The comparison principle for $\lcal(u)=0$ fails to hold in general: e.g. in the unit ball $\B^2$ in $\C^2$ the functions $u=|z_1|^2$, $v=1-|z_2|^2$ are different solutions of $\lcal(u)=0$ which coincide on ${\rm b}\B^2$. Neverthless, it holds for lower and upper solutions \cite[Prop. 4.3]{ST}. As we will point out below $\lp2.3\rp$, unicity for $\lp{\sf P}_g\rp$ is deeply related to certain hulls of the level sets $\{g={\rm const}\}$ (see below).
\er  
\subsection{Hulls and unicity}
We recall now the main results of \cite{ST} which were proved more generally for bounded {\rm P }-regular domains\footnote {A domain $\O$ in $\C^n$ is said to be  {\rm P }-{\em regular} if for every $z_0\in {\rm b}\O$ there is an open neighborhood $U$ of $z_0$ and  a continuous function $\rho$ in $U\cap\oli\O$ such that $\rho(z_0)=0$, $\rho(z)<0$ for $z\neq z_0$ and $\rho$ is p.s.h. in $U\cap\O$.} in $\C^2.$  
						
Let $\O$ be a bounded  $P$-regular domain in $\C^2$. If $u^\pm$ are the extremal solutions of the problem $({\sf P_g})$ then 
\beqnn\label{5(7)} 
\{u=0\}&\sbs&\{u^-\le c\}\cap \{u^+\ge c\}=\\
&&\Hat{\{g\le c\}}_{\pcal}\cap \Hat{\{g\ge c\}}_{\pcal}=\nonumber\\
&&\Hat{\{g=c\}}_{\pcal}\nonumber
\eeqnn
\cite[Corollary 5.3]{ST}.

As for unicity 
$$
\{u=c\}=\Hat{\{u=c\}}_{\pcal}=\Hat{\{g=c\}}_{\pcal}
$$
for every real $c\in[\min g,\max g]$  is a necessary condition (\cite[Cor. 5.4]{ST}).

A sufficient condition is that every hull $\Hat{\{g=c\}}_{\pcal}$, has empty interior (\cite[Prop. 5.5]{ST}).

Finally, the problem $({\sf P_g})$ has a unique solution if and only if the sets 
$\Hat{\{g=c_1\}}_{\pcal}$, $\Hat{\{g=c_2\}}_{\pcal}$ are mutually disjoint whenever $c_1\neqv c_2.$

 \br\label{r45}
 Observe that the levels $\{u=c\}$ of a weak solution of  $\lcal(u)=0$ might have non-empty interior: e.g. take on $\C^2$ the function $u$ defined by $u=0$ for ${\sf Re}z_1\le 0$ and $u={\sf Re}z_1$ otherwise. The existence of non-empty interior levels  solutions is intimately related with unicity (see \cite[Th. 5.6]{ST}). In particular, the Dirichlet problem on $\B^2$ considered in \ref{Exs} has a non-empty interior level solution.
 \er
 \bt\label{ETO}
Let $\O$ be a bounded domain in $\C^2$ and $u\in {\rm C}^0(\O)$ a solution of $\lcal(u)=0$ continuous up to the boundary. Let $g=u_{\vert{\rm b}\O}$ and assume that $S=\{g=0\}$ divides ${\rm b}\O$ into two connected components $S^\pm=\{g\gl 0\}$. Then the zero set of  
$u$ divides $\O$ into two connected components $\{u\gl 0\}$.  
\et
\demo
Denote $\oli U^\pm$, the connected components of $\{u\neq 0\}$ which contain $S^\pm$ respectively. Then $U^\pm={\oli U^\pm}\cap\O$ are the connected components of $\O\ssmi\{u=0\}$. By a contradiction let $V$ be a third connected component of $\O\ssmi\{u=0\}$. Since $u$ is vanishing on ${\rm b}V$, $\oli V$ cannot be compact otherwise, by the maximum principle 
$V\sbs\{u=0\}$. It follows that $\oli V\cap{\rm b}\O\sbs S$. Again, by the maximum principle $u$ must vanish on $V$: contradiction.
\enddemo
\section{The Main Theorem}\label{MT} 
Let $\O\sbs\C^2$ be a strongly pseudoconvex domain. Through this section $\oli M$ denotes a connected compact subset  of $\oli\O$ such that $S:=\oli M\cap{\rm b}\O$ is a connected subset dividing ${\rm b}\O$ into two connected components $S^{\pm}$.

We may assume that $S=\{g=0\}$ where $g$ is a continuous function $g:{\rm b}\O\to\R$ such that  $S^\pm=\{g \gl 0\}$  and $g$ has only two peak points $p\in S^+$, $q\in S^-$. Indeed, let $p\in S^+$, $q\in S^-$ respectively a  maximum and minimum point of $g$, $U^+\sbs S^+ $, $U^-\sbs S^- $ coordinate balls centered at $p$ and $q$ respectively with $\oli {U^\pm}\sbs S^\pm$. Ttake $\rho_1$, $\rho_2$ real smooth functions in ${\rm b}\O$ such that
\bit 
\item[1\rp] $\rho_1\le0$, $\rho_2\ge 0$, ${\rm supp}\rho_1\Sbs B_1$,
${\rm supp}\rho_2\Sbs B_2$;
\item[2\rp] $\min g>\rho_1(p)$, $\rho_1(z)>\rho_1(p)$ for $z\neqv p$ and $\max g<\rho_2(q)$, $\rho_(z)<\rho_2(q)$ for $z\neqv q$.
\eit
Then the function $g+\rho_1+\rho_2$ has the required properties.
 
 A function $g$ as above will be said a {\em minimal defining function}.

\bl\label{3l}
Let $\O\ssmi M$ have two connected components $\O^\pm$ and $\oli\ocal^\pm=\ocal(\oli{\O^\pm})$, $\hcal^\pm=\hcal(\O^\pm)$, $\pcal^\pm=\pcal(\O^\pm)$ (see Section \ref{Geo}. If $\acal$ denotes one of the symbols $\hcal$, $\pcal$,
$\oli\ocal$ then 
\bit
\item[\rm1\rp] $\Hat {S}_{\acal^\pm}\cap{\rm b}\O={\rm b}M$;
\eit
\bit                            
 \item[\rm2\rp] if $M:=\oli M\cap\O$ is l.m.p. the $\acal$-envelopes of ${\rm b}\O\cap{\rm b}\O^\pm$ and $\oli\O^\pm$ coincide.
\eit
\el
\demo
1) is a consequence of the fact that every point of  ${\rm b}\O$ is a peak point for $\hcal$. (This is clear if  $\O$ is strongly convex; the strongly pseuconvex case can be reduced to this one using an embedding theorem of Fornaess, see \cite[Theorem 8]{F}). 
\enddemo

\bl\label{F6}
Let $\O$ be a bounded domain in $\R^n$ with a ${\rm C}^1$ boundary. Let $K$ be a compact subset of $\oli\O$ such that $\Sigma=K\cap{\rm b}\O$ is a ${\rm C}^1$ connected $k$-dimensional submanifold of ${\rm b}\O$. Assume that $\Sigma$ is  the boundary of a $(k+1)$-dimensional submanifold $\G$ of ${\rm b}\O$. Then there exists a subset $S$ of $\oli\O$ with the properties
\bit
\item[i\rp] $S\cap{\rm b}\O=S\cap K=\Sigma$, $\Sigma\sbs\oli{S\ssmi\Sigma}$;
\item[ii\rp] $S\smi\Sigma$ is a  $(k+1)$-dimensional  ${\rm C}^1$ submanifold.
\eit    
\el
\demo
Let $d(\cdot,\cdot)$ be the euclidean distance in $\R^n$. Choose $\e$, $\d$ in such a way that the subsets 
$$
\big\{x\in\R^n:d(x,\Sigma)<\e\big\},\>\>\big\{x\in\R^n:d(x,{\rm b}\O)<\e\big\}
$$
define tubular neighborhoods of class ${\rm C}^1$ of $\Sigma$ and ${\rm b}\O$ respectively. Let $\{a_n\}_{n\in\N}$ be a strictly decreasing sequence of positive real numbers such that $a_n\to 0$ and $a_n<a<\e$. Let 
$$
U_n=\{x\in\R^n:d(x,\G)<a_n\}.
$$
We then choose a sequence $\{\e_n\}_{n\in\N}$ of positive real numbers such that 
$$
\e_n<\min\big\{\e, d(\oli U_n\ssmi U_{n+1}, K\big\}.
$$
It is immediately seen that there exists a ${\rm C}^1$ function $\psi:\R^+\to\R^+$ such that
\bit
\item[$\bullet$]  
$\>\>\>\>\psi<\d, \psi(x)$ {\rm is constant for} $x>a$ and $\psi(x)\to 0$ as $x\to 0 $
\item[$\bullet$] $\>\>\>\>\psi<\sum_{n\ge 0}{\e_n}{\bf 1}_{[a_{n+1},a_n]}$
\eit
where  ${\bf 1}_{[a_{n+1},a_n]}$ is the characteristic function of $[a_{n+1},a_n]$.
Let $x\in{\rm b}\O$ and $\nu(x)$ be the inner unit normal vector to ${\rm b}\O$ at $x$. We define the function $\rho:{\rm b}\O\to\oli\O$ by
$$
\rho(x)=x+\psi(d(x,\Sigma)\nu(x)).
$$
Clearly, $\rho$ is ${\rm C}^1$ and $S:=\rho(K)\cup\Sigma$ satisfies the conditions of the lemma $S\ssmi\Sigma$ is locally the graph of a ${\rm C}^1$- function).
\enddemo
 \bt\label{T10}
Let $\oli M$ be as above, $u\in{\rm C}^0(\oli\O)$ a solution of the Dirichlet problem $\lp{\sf P}_g\rp$ where $g$ is a minimal defining function for  ${\rm b}M$ . If $M$ is l.m.p. then
\beqn
a)\>\>\>\>\>&&M\sbseq\{u=0\}\sbseq\\
&&\{u^-\ge 0\}\cap\{u^+\le 0\}=\\
&&\Hat{\{g=0\}}_{\hcal}=\Hat{{\rm b}M}_{\hcal}.
\eeqn  
\nin  If $M$ is l.m.p. and ${\rm b}M$ is ${\rm C}^1$ then
\beqn
b)\>\>\>\>\>&&M\sbseq\{u=0\}\sbseq\Hat M_{\oli\ocal_1}\cup\Hat M_{\oli\ocal_2}=\\
&&\Hat {{\rm b}M}_{\oli\ocal_1}\cup\Hat {{\rm b}M}_{\oli\ocal_2}\sbs\Hat{{\rm b}M}_{\oli\ocal}.
\eeqn  
 \et
 \demo
 In view of \ref{5(7)}, to prove part a) it is enough to show that $M\sbs \{u=0\}$.
 \medskip
  
1) $\{u=g(p)\}=\{p\}$, $\{u=g(q)\}=\{q\}$.
\medskip 

\nin Indeed, $\{u=g(p)\}\cap{\rm b}\O=\{p\}$ since $g$ is minimal. Assume, by a contradiction, that  $\{u=g(p)\}\cap\O\neq\{p\}$. Then there exists a real plane $H$ in $\C^2\ssmi\{p\}$ such that for a connected component $\O_0$ of $\O\ssmi H$ we have 
$$
\O_0\cap\{u=g(p)\}\neqv\ES ,\>\>{\rm b}\O_0\cap\{u=g(p)\}=\ES.
$$
$\O_0\ssmi\{u=g(p)\}$ is Stein since $\O\ssmi\{u=g(p)\}$ is Stein (\ref{Exs}) . Moreover, by the extension theorem in \cite{LT}, every holomorphic function in $\O_0\ssmi\{u=g(p)\}$ extends through $\{u=g(p)\}$: contradiction. 
The proof for $\{u=g(q)\}=\{q\}$ is identical.\medskip

2) $M\sbs \{u=0\}.$

\medskip 
\nin If not there exists a point $z_0\in M$ such that $u(z_0)\neq0$, say $u(z_0)>0$. Since $u(p)=g(p)<0$ there exists $w\in \O_1$ such that $u(w)=0$.

Let 
$$
I=\big\{c\in[0,g(q)]:\{u=c\}\cap\O_1\neqv\ES\big\}.
$$
$0\in I$ and, by what proved in 1), $g(q)\notin I.$

 Let $c_0=\sup I$: $c_0>0$ and $c_0<g(q)$. If not we can find two sequences $\{c_\nu\}\sbs \R$, $\{\z_\nu\}\sbs\O_1$ such that $u(z_\nu)=c_\nu$, $c_\nu\to g(q)$, $\z_\nu\to\z_0\in\oli{\O_1}$. It follows that $u(\z_0)=g(q)$ and consequently, by 1) $\z_0=q$: contradiction since $\z_0\in\oli{\O_1}$. Therefore, $\{u=c_0\}\cap\O_1=\ES$ since $c_0=\sup I$ and $\{u=c_0\}\cap\oli\O_1\neq\ES$ otherwise $\{u=\d\}\cap\O_1=\ES$ for some $\d<c_0$ un contrast to  the hypothesis that  $c_0=\sup I$. Therefore, $\{u=c_0\}\cap\oli\O_1\sbs \oli M $ and since $u$ is vanishing on ${\rm b}M$, $\{u=c_0\}\cap\oli\O_1$ is a non empty compact subset of $M$: this violates   Kontinuit\"atsatz (see Section \ref{Geo}) since $\O\cap\{u>c_0\}$ is pseudoconvex  and $M$ is l.m.p..              
  The proof if $u(z_0)<0$ is similar. Thus $M\sbs\{u=0\}$.Taking into account \ref{5(7)} this proves part a) of the Main Theorem.
  
As for part b) we have to prove the second inclusion. Let $u(z)=0$. If $z\in\oli M$ then $z\in\Hat M_{\oli\ocal_1}\cup\Hat M_{\oli\ocal_2}.$ If $z\not\in\oli M$ then $z\in\O_1\cup\O_2$ since $u\neq0$ on ${\rm b}\O\smi\oli M$. Say $z\in\O_2$.  Assume, for a contradiction that $z\not\in M_{\oli\ocal_2}$ $_{\oli\ocal_2}$and apply Lemma \ref{F6} to $\oli\O_2$, $\G={\rm b}\O\cap{\rm b}\O_1$, $\Sigma={\rm b}M$ and $K=M_{\oli\ocal_2}\cup(\oli\O_2\cap\{u=0\})$. Let $\O_S\sbs\oli\O_2$ be the domain bounded by $S$ and $M$. In view of Lupacciolu's extension theorem \cite[Theorem 2]{L} every continuous CR function in $S\smi {\rm b}M$ extends holomorphically on $\O_ S$, an open subset containing $z_0$. In particular, this holds for every holomorphic function on the connected component $W$ containing $S$: contradiction since $W$ is pseudoconvex. This proves that $\{u=0\}\cap\oli\O:\sbs\Hat M_2$. In a similar way we show that $\{u=0\}\cap\oli\O_1\sbs\Hat M_1$. The proof of Theorem \ref{T10} is now complete. 
\enddemo
Given a connected compact subset $S$ of ${\rm b}\O$ let $\mcal(S)$ be the family of all subsets $\oli M$ of $\oli\O$ such that $M=\oli M\cap\O$ is  l.m.p. and $\oli M\cap{\rm b}\O=S$. $\mcal(S)$ is partially ordered by inclusion
\bt\label{C21} 
$\mcal(S)$ is non empty. If $S$ divides ${\rm b}\O$ into two connected components then
\bit
\item[\rm 1\rp] $\mcal(S)$ has a unique maximal element $M_S$; 
\item[\rm 2\rp] $M_S=\{u_g=0\}$ where $g$ is a minimal defining function for $S$ and  $u_g$ is a solution of the Dirichlet problem $\lp{\sf P}_g\rp$;
\item[\rm 3\rp] $M_S=\Hat S$. 
\eit
\et
\nin In particular, for every $\oli M\in\mcal(S)$ there is a continuous solution $u\neqv 0$ of the Levi equation such that $u=0$ on $M$. We will compare this result to \cite[Theorem 1.1]{FO}.

Theorem \ref{C21}  follows from
\bt\label{hulmp}
 Let $\O \Sbs\C^2$ be a strongly pseudoconvex bounded domain, $\hcal(\O)$ the Banach algebra $\ocal(\O)\cap {\rm C^0(\oli \O})$. For every non empty closed subset $S$ of ${\rm b}\O$ let 
 $$
 \Hat S =\{z\in{\oli \O}: \vert f(z)\vert\le\Vert f\Vert \forall f\in\hcal(\O)\}  
 $$
be the $\hcal(\O)-{\rm hull}$ of $S$. Then $\Hat S\ssmi S$ has the local maximum property.   
\et
\demo
1) Since $\O$ is strongly pseudoconvex, the {\rm Gelfand space} of $\hcal(\O)$ is homeomorphic to $\oli\O$.
Consider now the set of restrictions to $\Hat S$ of all the functions in $\hcal(\O)$.It is a subalgebra of ${\rm C}^0(\Hat S)$. Denote $\acal(\Hat S)$ the closure of this subalgebra with respect to the sup norm in ${\rm C}^0(\Hat S)$. The $\check{S}ilov$ boundary of $\acal(\Hat S)$ is contained in $S$ (actually equal, but we do not need this). The set $\Hat S$ is naturally contained in the Gelfand space  of $\acal(\Hat S)$ via points evaluations.

2) The Gelfand space of $\acal(\Hat S)$ is equal to $\Hat S$.

\nin Let $\chi:\acal(\Hat S)\to\C$ be a character of $\acal(\Hat S)$. Let $r:\hcal(\O)\to \acal(\Hat S)$  be the restriction homomorphism sending $f\in\hcal(\O)$ to $f_{|\Hat S}$. The composition $\chi\circ r$ is a character of  $\hcal(\O)$ and so, by 1), there is a unique point $(z_1,z_2)\in\oli\O$ such that for every  $f\in\hcal(\O)$, $\chi(f_{|\Hat S})=f(z_1,z_2).$ Since $\chi$ has norm $1$, we get
$$
\vert f(z_1,z_2)\vert=\vert\chi(f_{|\Hat S})\vert\le\sup \vert f_{|\Hat S}\vert=\sup\vert f_{|S}\vert
$$
for every $f\in\hcal(\O)$ and this shows that $(z_1,z_2)$ belongs to $\Hat S$.

In order to ends the proof we observe that by Rossi's local maximum modulus principle, if $f\in\hcal(\O)$, $\vert f\vert$ cannot have local maximum at any point of $\Hat S\ssmi S$. Suppose, by contradiction that $\Hat S\ssmi S$ is not l.m.p. Then there exists a holomorphic polynomial $Q(z_1,z_2)$ such that $Re Q _{|\Hat S\ssmi S}$ has local maximum at some point of $\Hat S\ssmi S$ ( Section \ref{INTR}, LM1). Then the modulus of $f=\exp Q$ would have a global maximum at some point of $\Hat S\ssmi S$ and $f$ is in $\hcal(\O)$; contradiction.
\enddemo
  
\bexc
Let $\mathbb T^2$ be the torus $\{|z_1|=|z_2|=1\}$ in the sphere $\mathbb S^3=\{|z_1|^2+|z_2|^2=2\}$. Then the solid torus $\widetilde{\mathbb T}^2=\{\vert z_1\vert\le 1,\vert z_2\vert\le1\}$ is maximal in $\mcal(\mathbb T^2).$
\eexc

\bc \label{T11}
Assume that $\Hat {{\rm b}M}_ {\pcal(\oli\O)}=M$. Then there exists a base of neighborhoods $\{U_n\}_{n\in\N}$ of $M$ for the relative topology of $\O$ such that $\Til U_n:=U_n\cap\O$, $n\in\N$, are Stein domains.
In particular, $\Int{M}=M\ssmi{\rm b}M$ is union of an increasing sequence of Stein compacts. 
\ec
\bex Let $\mathbb T$ be the torus $\{|z_1|=|z_2|=1\}$. Then $\{|z_1|\leq 1, |z_2|=1\}$ and $\{|z_2|\leq 1, |z_1|=1\}$ are Levi flat hypersurfaces with boundary $\mathbb T$ which are also Stein compact sets. On the other hand, $X=\{|z_1|=|z_2|,|z_1|\leq 1\}$ is a (singular) Levi flat hypersurface with boundary $\mathbb T$, but it is not a Stein compact set.
\eex

\demo
Define $V_\epsilon=\{||z_1|-|z_2||<\epsilon, |z_1|,|z_2|\leq 1\}$. The sets $V_\epsilon$ are a fundamental system of neighborhoods of $X$. It is enough to show that, for all $\epsilon>0$, the holomorphic completion of $V_\epsilon$ includes the bidisc $\{|z_1|\leq 1, |z_2|\leq 1\}$. For $M\leq 1$ put $A_M=V_\epsilon\cup\{|z_1|,|z_2|<M\}$, and let $M'$ be the supremum of the values of $M$ such that $A_M$ is contained in the holomorphic completion of $V_\epsilon$. Note that $M'\geq \epsilon/2$ since $A_{\epsilon/2}=V_\epsilon$. We need to show that $M'=1$: assume by contradiction that $M'<1$.

Consider the Hartogs set $H=H_1\cup H_2$, where
\[H_1=\{z_1=\alpha, |z_2|\leq M'+\epsilon/12 \},\]
\[H_2= \{|z_1-\alpha|\leq \epsilon/3, |z_2|= M'+\epsilon/12\} \]
with $|\alpha|=M'-\epsilon/4$. A holomorphic function defined on a neighborhood of $H$ extends to the bidisc
\[P=\{|z_1-\alpha|< \epsilon/3, |z_2|<  M'+\epsilon/12\}. \]
\begin{center}
   \includegraphics[width=8cm]{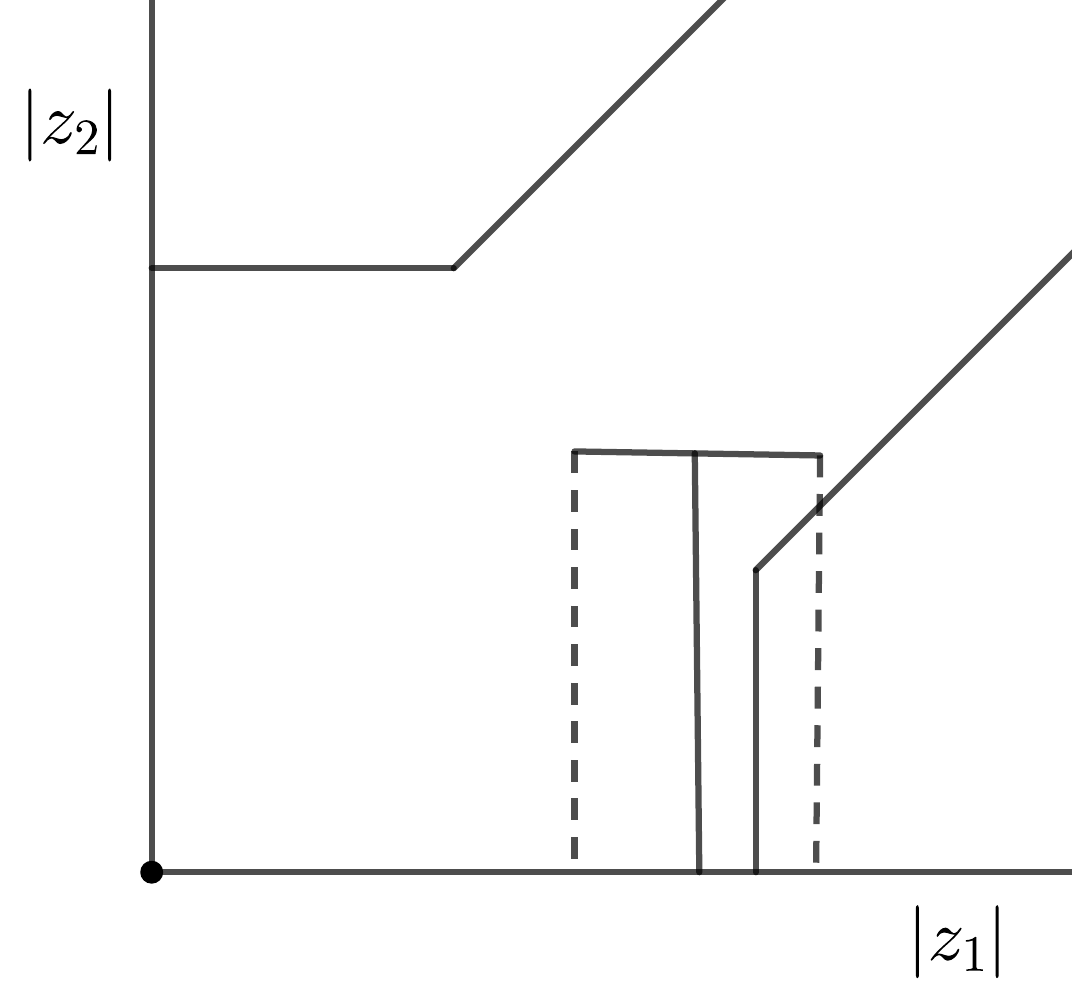} 
\end{center}
On the other hand, $H\subset A_{M'}$. Indeed, for all $(z_1,z_2)\in H_1$ we have either $|z_2|<M'$ or $||z_1|-|z_2||\leq\epsilon/3$, while for $(z_1,z_2)\in H_1$ we have $||z_1|-|z_2||\leq 2\epsilon/3$.
Thus, the holomorphic completion of $V_\epsilon$ contains $A_{M'}\cup P$, and taking all $\alpha$ with $|\alpha|=M'-\epsilon/4$, it follows that it contains 
\[A_{M'}\cup \{|z_1|<M'+\epsilon/12,|z_2|\leq |z_1| \}\]
and finally, exchanging the role of $z_1$ and $z_2$, we find that it contains $A_{M'+\epsilon/12}$, against the definition of $M'$.
\enddemo

\section{Geometric properties of Levi convex graphs}
\subsection{The Levi operator for graphs}\label{LG}
Through this subsection $D$ is a bounded domain in $\C_{z_1}\tms\R_{x_3}$. For every function $v:D\to\R_{x_4}$ we set
\begin{center} 
$
\G(v)=\{x_4=v\}, \G(v)^\pm=\{x_4\gl v\}.
$
\end{center}
The Levi condition \ref{5(6)} for $\G(v)$ writes
\beqnn\label{L3}%numerata
&&\lsf(v)=(1+v_3^2)(v_{11}+v_{22})+(v_1^2+v_2^2)v_{33}+\\
&&2(v_2-v_1u_3)v_{13}-2(v_1+v_2v_3)v_{23}=0\nonumber
\eeqnn
(where $v_i:=\p v/\p x_i$, $v_{ij}:=\p^2 v/\p x_i\p x_j$). 
\nin Even $\lsf(v)$ is proper and elliptic degenerate.  

More generally, we consider the differential operator
\beqnn\label{L8}%numerata
&&\lsf(v;k)=(1+v_3^2))(v_{11}+v_{22})+(v_1^2+v_2^2)v_{33}+\\
&&2(v_2-v_1v_3)v_{13}-2(v_1+v_2v_3)v_{23}+\nonumber\\
&&k(\cdot;v)(1+\vert Dv\vert^2)^{3/2}\nonumber
\eeqnn
where $k=k(x;t)$ is continuous in $D\tms\R_{x_4}$.
 Let $v\in{\rm C}^0(\oli D)$. Then
 \bit
\item[a\rp] we say that $v$ is a {\em weak subsolution} of $\lsf(v;k)=0$ if for every $y\in D$ and $\phi\in{\rm C}^\IN(D)$ such that $v-\phi$ has a local maximum at $y$ one has $\lsf(\phi;k)(y)\ge 0$;
\item[b\rp]we say that $v$ is a {\em weak supersolution} of $\lsf(\phi;k)=0$ if for every $y\in D $ and $\phi\in{\rm C}^\IN(D)$ such that $v-\phi$ has a local minimum at $y$ one has $\lsf(\phi;k)(y)\le 0$
\eit
$v$ is said to be a {\em weak solution} of $\lsf(v;k)=0$ if it is both a weak subsolution and a weak supersolution.

\br\label{ERM}
The Dirichlet problem was considered, more generally, for the ``complete Levi operator''
$${\lcal
}(u;k)=\sum_{\a,\b=1}^2(\delta_{\a\b}{\vert\partial
u\vert^{2}-u_{\bar\a}u_{\b}})u_{\a\bar\b}-k(\cdot,u)\vert\p u\vert^3. 
$$
where $k;\O\tms\R\to\R$. Under suitable conditions for $k$, if $g\in{\rm C}^{2,\a}(\rm b\O)$ the problem has a solution $u\in{\rm Lip}(\oli\O)$  (\cite[Th. 2.3]{ST}).
\er

The following identities are immediate:
\bit
\item[6\rp] $\lsf(v)=-\lcal(x_4-v)=\lcal(v-x_4)$;	
\item[7\rp]  $\lsf(-v)=\lcal(-v-x_4)=-\lcal(x_4+v)$;
\item[8\rp] $\lcal(v-x_4)=\lsf(v;\til k)$ where\\ $\til k(x,t)=k(x,v-x_4)$.
\eit
Moreover, if $v$ is a continuos subsolution (supersolution, solution) of $\lsf(v)=0$ in $D$, then $v-x_4$ is a subsolution (supersolution, solution) of $\lcal(u)=0$ on $D\tms\R_{x_4}$. From 1), 2), 3) and (\cite[Corollary 3.2]{ST}) we get the following:
\bit
\item[9\rp] if $k\ge 0$ and $\lsf(v;k)\le0$, then $\G(v)^-$ is pseudoconvex;	
\item[10\rp] if $k\le 0$ and $\lsf(v;k)\ge 0$ then $\G(v)^+$ is pseudoconvex. 
\eit 
The graph $\G(v)$ of a continuous function $v:D\to\R_{x_4}$ is called a {\em family of analytic discs} if it is a disjoint union of analytic discs.
\bp\label{P}
A graph  $\G(v)$ of a continuous function $v:D\to\R_{x_4}$ is a family of analytic discs if and only if $\lsf(v)=0.$
\ep
\demo
Assume that $\G(v)$ is a family of analytic discs. Then $\G(v)$ is l.m.p. . If not, by conditin LM2 there exist $z^0\in X$, $\e>0$ a ball $B(z^0,r)$ and a strictly plurisubharmonic function $\psi$ on $B(z^0,r)$ such that $\psi(z^0)=0$ and $\psi(z)\le-\e\vert z-z^0\vert^2$ for $z\in X\cap B(z^0,r)$. Then, on ananalytic disc $D\sbs X\cap B(z^0,r)$ through $z^0$ the plurisubharmonic function $\psi_{\vert D}$ violates the maximum principle. Therefore $\G(v)$ is l.m.p. and consequently all the hypersurfaces $u=x_4-v=const$ are l.m.p. From \cite[Theorem 3.3] {ST} it follows that $\lcal(u)=0$ and consequently that $\lsf(v)=0.$

Conversely, let $\lsf(v)=0.$. Then $\lcal(u)=0$ and, again by \cite[Theorem 3.3] {ST}  $D\tms\R_{x_4}\ssmi\G(v)$ is locally pseudoconvex. By the Main Theorem of \cite{PS} $\G(v)$ is a family of analytic discs. 
\enddemo
 For any compact subset $K$ of $\C^2$ we denote $\Hat K$ its polynomial envelope.   If $C$ is a subset containing $K$ we denote ${\rm hull}_{\ocal(C)}K$ the envelope of $K$ with respect to the algebra $\ocal(C)$.
\bp\label{P1}
Let $D$ be bounded with ${\rm b}D$ of class $C^{2,\a}$, $0<\a<1$  and $Dì\tms\R_{x_4}$ strictly pseudoconvex. Then 
\bit
\item[\rm1\rp] For every $g\in C^0({\rm b}D)$ there exists a unique $v\in C^0(\oli D)$ such that $\lsf(v)=0$ in $D$ and $v=g$ on ${\rm b}D$;  
\item[\rm2\rp] $\G(v)={\rm hull}_{\ocal(\oli D\tms\R_{x_4})} \G(g)$;
\item[\rm3\rp]if $\oli D\tms[l-,l]$ is polynomially convex for some $l>\max\limits_{\bar D}v$, in particular if $\oli D$ is convex, $\G(v)=\Hat {\G(g)}$;
\item[\rm4\rp] $\G(v)\ssmi \G(g)$ is a family of analytic discs.
\eit
\ep

\nin 1), 2) were proved in \cite{ST2} (Th. 3.1 and Prop. 7).  

\nin 3) and 4) were proved by Shcherbina (\cite[Main Theorem]{Sh}.4 ) for strictly convex domains $D\tms\R$ and by Chirka and Shcherbina (\cite[Th. 2]{CS}) when $D\tms\R$ is strictly pseudoconvex. 
\br\label{R8}
In the proof of \cite[Main Theorem]{Sh}.4 )  is crucial the existence of an analytic foliation on $\G(v)$ when $g$ is smooth. This is one of the main results of Bedford and Klingenberg in \cite{BK}{\lp} cfr. Theorem 3{\rp}. An analogous result was proved in \cite{CT} for the solutions of the Levi equation for almost complex structures {\lp} {\cite[Th. 1]{CT}} for $f=0${\rp}.
\er

 \subsection{Min-max principle}\label{M-m}
 \bp
Let $\O\sbs\R_x$ be a bounded pseudoconvex domain, $v\in{\rm C}^0(\oli\O)$ a subsolution of $\lsf(v;k)=0$ in $\O$. Assume that $k\le 0$. Then
\bit 
\item[\rm1)] $\max\limits_{\oli\O}\,v=\max\limits_{{\rm b}\O}\,v;$  
\item[\rm 2)] if 
$$
M:=\G(v)\cap\{x_4=\max\limits_{\oli\O}\,v\}\neqv\ES
$$
then for every complex line $l\sbs\{x_4=\max\limits_{\oli\O}\,v\}$ the subset $l\cap(M\ssmi{\rm b}\G(v))$ has no compact connected component.    
\eit
\ep
\demo
1) Since $k\le 0$ $v$ is a subsolution of $\lsf(v)=0$, therefore $u=v-x_4$ is a subsolution of $\lcal(u)=0$ in $Q_c:=\O\tms(-c,c)$, continuous in $\oli{Q_c}$. 
Assume, by contradiction, 
$$
m:=\max\limits_{\oli\O}\,v>\mu:=\max\limits_{{\rm b}\O}\,v.
$$
Then, for some $\a$,
 $$
 \{x_4=\a\}\cap\G(v)^-=\ES
 $$ 
 and 
 $$
 K:=\{x_4=\a\}\cap\G(v)
 $$
 is compact in $\{x_4=\a\}\cap Q_c$. Let $V$ a relatively open neighborhood of $K$ in $\{x_4=\a\}$ with $\oli V$ compact. For every $\nu\in\N$, $\nu\gg 0$, let 
 $$
 X_\nu=V\tms\{\a<x_4<\a+1/\nu\}.
 $$
 All the $X_\nu$ are l.m.p, the conditions a), b) of Section \ref{Geo} are fulfilled  but $\G(v)^+$ is pseudoconvex: contradiction.

2) If $l\cap(M\ssmi{\rm b}\G(v))$ had a non empty compact connected component $K'$, it would exist a relative open neighborhod $D\sbs l$ of $K'$ in $l$ such that ${\rm b}D\sbs\G(v)^+$. Moving $D$  we would get then a family $\{X_\nu\}$ as above which violates again \cite[Corollary 3.2]{ST}.     
\enddemo
The following min-max principle holds true for (weak) subsolutions and supersolutions
\bt\label{L4}
Let $u\in{\rm C}^0(\oli\O)$ and $k=k(x)$. If $u$ is a subsolution (supersolution) of $\lsf(k;v)=0$ and $k<0$ ($k\ge0$) then $v$ has no local maximum (minimum) in $\O$. If  $v$ is a subsolution (supersolution) of $\lsf(k;v)=0$ and $k\le 0$ ($k\ge0$) then  
$$
\max\limits_{\oli\O}v=\max\limits_{{\rm b}\O}v\>\> (\min\limits_{\oli\O}v=\min\limits_{{\rm b}\O}v).
$$
In particular if $k=0$ and $v$ is a solution then for every $x\in\oli\O$
$$
\vert v(x)\vert\le\max\limits_{{\rm b}\O}\vert v\vert
$$
\et

\subsection{Dirichlet problem} In \cite[Theoren 4]{ST2} the following is proved. Let $\O$ be a bounded domain, ${\rm b}\O\in C^{2,\a}$, $g\in{\rm }C^{2,\a}({\rm b}\O)$, $0<\a<1$ and $k\in{\rm }C^1(\oli\O\tms\R)$. Under suitable conditions on $k$ the Dirichlet problem
\[
\hspace{-1cm}({\sf P}^\ast_{k,g})\hspace{1cm}\begin{cases}
\lsf(v;k)=0 \>\> & {\rm in}\>\> \O\\
v=g\>\> & {\rm on}\>\> {\rm b}\O
\end{cases}
\]
has a weak solution $u\in{\rm Lip}(\oli\O)$. In particular, the problem $({\sf P}^\ast_{0,g})$ has a ${\rm Lip}(\oli \O)$ solution.
\subsection{Geometric properties of solutions}. For $k=0$ the solutions of $({\sf P}^\ast_{0,g})$ give particular solutions of $({\sf P}_g)$ by $v\mapsto x_4-v$ so in this case  we are reduced to Section \ref{Geo}.
\bp\label{p12}
Let $\O\sbs\R_x$ be a bounded domain, $v\in{\rm C}^0(\O)$ a solution of $\lsf(v;k)=0$ in $\O$ such that the hypersurface $\G(v)$ is complete. Then 
\bit 
\item[\rm 1)] if $k\ge 0$, $\inf\limits_\O v=-\IN$ and $\sup\limits_\O v<\IN$
\item[\rm 2)] if $k\le 0$, $\inf\limits_\O v>-\IN$ and $\sup\limits_\O v=+\IN$
\eit
In particular, no complete l.m.p. continuous graph exists on $\O$.
\ep
\demo
Since $\G(v)$ is complete, $Z:=\{x\in\O:v(x)=0\}$ is compact, hence  there exists a neighborhood $U$ of ${\rm b}\O$ such that
$v$ has a constant sign on $U\cap{\rm b}\O$, hence either
\bit 
\item[\rm1')] if $v\ge 0$ then $\inf\limits_\O v<-\IN$ and $\sup\limits_\O v=\IN$
\eit
or
\bit
\item[\rm 2')] if $v\le 0$ then $\inf\limits_\O v=-\IN$ and $\sup\limits_\O v<+\IN.$
\eit
In the former (latter) $v$ has a (maximum) minimum value so that for some $c\in\R$ the hyperplane $x_4=c$ contains a bounded domain $X$ such that ${\rm b}X$ lies in $\G^-(v)$ ($\G^+(v)$) but $X$ not. By Kontinuit\"atsatz (see \ref{Geo}) the domain $\G^-(v)$ ($\G^+(v)$) is not pseudoconvex. Since $u=x_4-v$ is a solution of $\lcal(u)=0$, Theorem \ref{T2} shows that the case $k\equiv 0$ is not allowed.

If $k\ge 0$, from $\lsf(v;k)=0$ we get $\lsf(v)\le 0$ and consequently $\lcal(x_4-v)=-\lsf(v)\ge 0$ i.e $x_4-v$ is a subsolution of  $\lcal(u)=0$. By \cite[Corollary 3.2]{ST} $\G^-(v)$ is pseudoconvex therefore, by what is preceding,  
$$
\inf\limits_\O v=-\IN\>{\rm and}\> \sup\limits_\O v<+\IN.
$$ 
This proves 1). The proof of 2) is similar. 

If $\G(v)$ is an l.m.p. continuous graph, then $v$ is a solution of $\lsf(v)=0$ (\ref {T2}) therefore 1) and 2) show that $\G(v)$ cannot be complete.    
\enddemo

\bt\label{T12} Let $B=B(R)\sbs\R^3$ the open ball of radius $R$ centered at the origin and  $u\in{\rm Lip}_{\rm loc} (B)$ be a solution of $\lsf(u;k)=0$ with $k=k(x)$. Then
\bit 
\item[\rm i)] if $k\ge 0$, $\inf\limits_Bk\le1/R$;
\item[\rm ii)] if $k\le 0$, $\sup\limits_B k>-1/R$ and $\sup\limits_\O v=+\IN$
\eit
In particular if $R=+\IN$ 
and $k\geq 0$ ($k\leq 0$), then $k$ cannot have a positive (negative) infimum (supremum).
\et
\demo
Let $k\le0$. It is sufficient to prove the inequality for every $R'<R$ therefore we may assume $u\in{\rm Lip}(\oli B)$. We may also assume $\min\limits_{{\rm b}B}(u-v)=0$. 

\nin Let $v=(R^2-\vert x\vert ^2)^{1/2}$: $v$ is a solution of $\lsf(v;1/R)=0$. The subset $\{u\le v\}\cap B$ is non empty. Indeed, consider a point $x_0\in{\rm b}B$ such that $u(x_0)=0$ and $\bf r$ the normal straight line to ${\rm b}B$ at $x_0$. Since $u\in{\rm Lip}(\oli B)$ and $u(x_0)=v(x_0)=0$, for every $x\in {\bf r}\cap B$ such that $u(x)>v(x)$ we have 

$$
C\ge(R-\vert x\vert)^{-1}u(x)\ge R-\vert x\vert)^{-1}v(x)=
$$
$$
((R+\vert x\vert)^{1/2}(R-\vert x\vert)^{-1/2}.
$$
where $C$ is a Lipschitz constant. Since $u\in{\rm Lip}(\oli B)$, this shows that in ${\bf r}\cap B$, near $x_0$, $u(x)<v(x)$, in particular that $u-v$ has minimum point on $B$, $u-v$ being non negative on ${\rm b}B$. Let $y\in B$ a local minimum point  for $u-v$. Then since  $v$ is smooth and $u$ is in particular a supersolution of $\lsf(u;k)=0$ we must have $\lsf(v)(y)+k(y)(1+\vert Dv\vert^2)(y)^{3/2}\le0$. On the other hand, since $\lsf(v;R^{-1})=0$ we have 
$\lsf(v)(y)+R^{-1}(1+\vert Dv\vert^2)(y)^{3/2}=0$ whence  $k(y)\le R^{-1}$.

The case $k\le 0$ reduces to the previous one observing that the function $\Tilde u:=-u(x_1,x_2,-x_3)$ is a ${\rm Lip}_{\rm loc} (B)$ solution of $\lsf(\Til u;-k)=0$.
\enddemo
\bc\label{T38}
Let $B=B(+\IN)=\R^3$. Then
\bit
\item[\rm iii)] if $\inf\limits_{\R^3}k>0$ or $\sup\limits_{\R^3} k<0$ the equation $\lsf(u;k)=0$ has no 
${\rm Lip}_{\rm loc}$ solution;
\item[\rm iv)] the equation $\lsf(u;k)=0$ with $k$ constant has ${\rm Lip}_{\rm loc} (\R^3)$ solutions if and only if $k\equiv 0$.
\eit   
\ec

\end{document}